\documentclass[12pt]{article}
\usepackage{amsmath,amsthm}
\usepackage{amssymb,latexsym}
\usepackage{enumerate}

\def\Q{{\mathbb Q}}
\def\Z{{\mathbb Z}}

\newtheorem{lemma}{Lemma}
\newtheorem{theorem}[lemma]{Theorem}

\title{
Calculating all elements of minimal index \\
in the infinite parametric family of simplest quartic fields
}
\author{
Istv\'{a}n Ga\'{a}l (Debrecen) \thanks{
        Research supported in part by K75566 and K100339 from the
        Hungarian National Foundation for Scientific Research
        and
        by the TAMOP 4.2.1./B-09/1/KONV-2010-0007 
        project implemented through the New Hungary Development Plan co-financed 
        by the European Social Fund, and the European Regional Development Fund
       }\;
and G\'abor Petr\'anyi (Debrecen)
\\
}

\begin{document}

\maketitle
\thispagestyle{empty}

\noindent
Mathematics Subject Classification: Primary 11Y50; Secondary 11D25,11R04\\
Key words and phrases: simplest quartic fields, power integral basis, monogenity\\

\begin{abstract}
\noindent 
It is a classical problem in algebraic number theory to decide if a
number field is monogeneous, that is if it admits power integral bases. 
It is especially interesting to consider this question in an infinite
parametric familiy of number fields. In this paper we consider the
infinite parametric family of simplest quartic fields $K$ generated by a root
$\xi$ of the polynomial $P_t(x)=x^4-tx^3-6x^2+tx+1$, assuming that $t>0$,
$t\neq 3$ and $t^2+16$ has no odd square factors. In addition to generators
of power integral bases we also calculate
the minimal index and all elements of minimal index in all fields in this 
family.
\end{abstract}

\section{Introduction}

There is an extensive literature (cf. \cite{gaalbook}) 
of {\it monogenic fields} $K$, that is algebraic
number fields of degree $n$ having 
a {\it power integral basis} $1,\vartheta,\ldots,\vartheta^{n-1}$.
This is the case exactly if the {\it index} of $\vartheta$, that is 
\[
I(\vartheta)=(\Z_K^+:\Z[\vartheta]^+)
\]
(where $\Z_K$ denotes the ring of integers of $K$) is equal to 1.

The first author developed algorithms for calculating generators of power integral bases
(cf. \cite{gaalbook}) and also succeeded to determine all possible generators of
power integral bases in some infinite parametric families of number fields,
see I.Ga\'al and M.Pohst \cite{gpc5}, I.Ga\'al and G.Lettl \cite{gaallettl}.

Among the best known infinite parametric families of number fields are the 
family of simplest cubic (see \cite{shanks}), simplest quartic and simplest 
sextic fields (see \cite{lpv}). These are exactly the families of totally real cyclic 
number fields, having a transformation of type $z\rightarrow \frac{az+b}{cz+d}$ 
as the generator of the Galois group.

In the present paper we deal with the family of {\em simplest quartic fields}, that is
$K=\Q(\xi)$, where $\xi$ is a root of the polynomial
\begin{equation}
P_t(x)=x^4-tx^3-6x^2+tx+1
\label{pp}
\end{equation}
where $t\in\Z,t\neq 0,\pm 3$. 
This family was considered by M.N.Gras \cite{gras}, A.J.Lazarus \cite{lazarus} and many
other authors.
Using the integral basis constructed by H.K.Kim and J.H.Lee \cite{kim},
P.Olajos \cite{olaj} showed that $K$ has power integral bases only in two exceptional cases.
He used the method of I.Ga\'{a}l, A.Peth\H{o} and M.Pohst \cite{gppsim} to solve the
index form equations. 

In the present paper we determine the minimal indices and all elements of minimal index
in the fields belonging to the infinite parametric family of simplest quartic fields.
Our basic tool is the method \cite{gppsim}, involving extensive formal 
calculations using Maple \cite{maple} and the resolution of a great number of
special Thue equations using Kash \cite{kash}. 

Note that B.Jadrievi\'c \cite{jadrievic}, \cite{jadrievic2}
has determined the minimal indices and all elements of minimal index
formerly  in certain families of bicyclic biquadratic fields. 
In those fields the index form equation splits into the product of
three quadratic forms, which makes the problem much easier.

\section{Simplest quartic fields}

For simplicity let $t>0,t\neq 3$ and $\xi$ a root of the polynomial (\ref{pp}).
We also assume that $t^2+16$ is not divisible by an odd square since this
was needed by H.K.Kim and J.H.Lee \cite{kim} to determine the integral bases.
(M.N.Gras \cite{gras} showed that $t^2+16$ represents infinitely many
square free integers. This implies that there are infinitely many 
parameters $t$ with the above properties.)

\begin{lemma} (H.K.Kim and J.H.Lee \cite{kim})\\
Under the above assumptions an integral basis of $K$ is given by
\[
\begin{array}{rl}
\left(1,\xi,\xi^2,\frac{1+\xi^3}{2}\right) & {\rm if}\; v_2(t)=0,\\
\left(1,\xi,\frac{1+\xi^2}{2},\frac{\xi+\xi^3}{2}\right) &  {\rm if}\; v_2(t)=1,\\
\left(1,\xi,\frac{1+\xi^2}{2},\frac{1+\xi+\xi^2+\xi^3}{4}\right) &  {\rm if}\; v_2(t)=2,\\
\left(1,\xi,\frac{1+2\xi-\xi^2}{4},\frac{1+\xi+\xi^2+\xi^3}{4}\right) &  {\rm if}\; v_2(t)\geq 3.\\
\end{array}
\]
\label{kimlemma}
\end{lemma}

P.Olajos \cite{olaj} determined all generators of power integral bases (up to translation
by elements of $\Z$).

\begin{lemma} (P.Olajos \cite{olaj})\\
Under the above assumptions power integral bases exist only for $t=2$ and $t=4$.
All generators of power integral bases are given by

\begin{itemize}

\item
$t = 2$,\;\; $\alpha=x\xi+y\frac{1+\xi^2}{2}+z\frac{\xi+\xi^3}{2}$
where\\
$(x, y, z) = (4, 2,-1), (-13,-9, 4), (-2, 1, 0),
(1, 1, 0),\\ (-8,-3, 2), (-12,-4, 3), (0,-4, 1),
(6, 5,-2), (-1, 1, 0), (0, 1, 0)$

\item
$t = 4$,\;\; $\alpha=x\xi+y\frac{1+\xi^2}{2}+z\frac{1+\xi+\xi^2+\xi^3}{4}$
where\\
$(x, y, z) = (3, 2,-1), (-2,-2, 1),
(4, 8,-3), (-6,-7, 3), \\(0, 3,-1), (1, 3,-1)$.

\end{itemize}
\label{olajlemma}
\end{lemma}

\section{Elements of minimal index in the family of simplest quartic fields}

If $K$ admits a power integral basis, then its generator has index 1.
Otherwise we call $m$ the {\it minimal index} of $K$ if $m$ is the least
positive integer such that $\alpha\in\Z_K$ exists with 
\[
I(\alpha)=m.
\]
If there are no power integral bases then it is important to determine
the minimal index of the field $K$ and all elements of minimal index.
It is easily seen from Lemma \ref{kimlemma} that $\xi$ has index 2,4,8,16
according as $v_2(t)=0,1,2,\geq 3$, respectively. We shall see that in some cases
there are elements of smaller index, as well. Moreover we determine all elements
of minimal index.

In the following theorem the coordinates of the elements are given in
the integral bases of Lemma \ref{kimlemma}.
We display only the last three coordinates of the elements (and omit the first one)
since the index is translation invariant.

Our main result is the following theorem. 

\begin{theorem}\mbox{}\\
\label{ththth}
\noindent
Assume that $t>0$, $t\neq 3$ and $t^2+16$ has no odd square factors.
Except the parameters $t=2,4,8,12,16,20,24,28,32$, for different values of
$v_2(t)$ the minimal indices $m$ of the field $K$ and all 
elements of minimal index are listed below: 

\[
\begin{array}{|l|}
\hline
v_2(t)=0, \;\;  m=2\\ 
(1,0,0),(6,t,-2),
\displaystyle{\left(\frac{5+t}{2},\frac{-1+t}{2},-1\right),
\left(\frac{5-t}{2},\frac{1+t}{2},-1\right) }         \\ \hline
v_2(t)=1, \;\;  m=4\\ 
(1,0,0),(7,2t,-2),
\displaystyle{\left(\frac{6+t}{2},-1+t,-1\right),
\left(\frac{6-t}{2},1+t,-1\right)}      \\ \hline
v_2(t)=2, \;\;  m=8\\ 
(1,0,0),(7,2+2t,-4),
\displaystyle{\left(\frac{6+t}{2},t,-2\right),
\left(\frac{6-t}{2},2+t,-2\right)}     \\ \hline
v_2(t)\ge 3,\;\; m=16  \\ 
(1,0,0),(9+2t,4-4t,-4),
\displaystyle{\left(\frac{10+t}{2},-4-2t,-2\right),
\left(\frac{6+3t}{2},-2t,-2\right)}        \\ \hline
\end{array}
\]

For $t=2,4,8,12,16,20,24,28$ the minimal indices and
all elements of minimal index are listed below: 

\[
\begin{array}{|l|}
\hline
t=2, \;\;  m=1\\ 
(4, 2,-1), (-13,-9, 4), (-2, 1, 0),(1, 1, 0),(-8,-3, 2), (-12,-4, 3),\\
(0,-4, 1), (6, 5,-2), (-1, 1, 0), (0, 1, 0)\\ \hline
t=4, \;\;  m=1\\ 
(3, 2,-1), (-2,-2, 1), (4, 8,-3),(-6,-7, 3), (0, 3,-1), (1, 3,-1)\\ \hline
t=8, \;\;  m=4\\ 
(-7,8,1),(17,-28,-3),(-8,8,1),(5,-10,-1),(20,-26,-3),(-4,10,1)\\ \hline
t=12, \;\;  m=3\\ 
(10,6,-1),(8,6,-1),(-4,-20,3),(-2,7,-1),(16,19,-3),(4,-7,1)\\ \hline
t=16, \;\;  m=8\\ 
(-14,16,1),(27,-52,-3),(-13,16,1),(6,-18,-1),(-34,50,3),(-7,18,1)\\ \hline
t=20, \;\;  m=5\\ 
(14,10,-1),(0,-32,3),(12,10,-1),(6,-11,1),(-8,11,-1),(-20,-31,3)\\ \hline
t=24, \;\;  m=12\\ 
(-19,24,1),(37,-76,-3),(-20,24,1),(-9,26,1),(8,-26,-1),(-48,74,3)\\ \hline
t=28, \;\;  m=7\\ 
(-18,-14,1),(-16,-14,1),(-4,44,-3),(-24,-43,3),(-13,14,1),\\(10,-15,1)\\ \hline
\end{array}
\]

\vspace{0.5cm}
\noindent
For $t=32$ in addition to those given for $v_2(t)\geq 3$ 
there are further elements of index 16 having coefficients 
\[
\begin{array}{|l|}
\hline
t=32, \;\;  m=16\\ 
(-26,32,1),(-25,32,1),(47,-100,-3),(11,-34,-1),(-10,34,1),\\(62,-98,-3)\\ \hline
\end{array}
\]

\label{mainth}
\end{theorem}

\section{Proof of Theorem \ref{ththth} }

In this section we list the results that we need to prove Theorem \ref{ththth} and
then we describe its proof.

\subsection{The corresponding family of Thue equations}

In our calculation we shall use the result giving all solutions $p,q\in\Z$
of the infinite parametric family of Thue equations
\begin{equation}
F_t(p,q)=p^4-tp^3q-6p^2q^2+tpq^3+q^4= w
\label{thuefamily}
\end{equation}
for given $w\in\Z$.

These equations were considered by G.Lettl and A.Peth\H{o} \cite{lp}
and by G.Lettl, A.Peth\H{o}, P.Voutier \cite{lpv}. 
Note that if $(p,q)$ is a solution, then so also is $(-p,-q)$
but we list only one of them. G.Lettl and A.Peth\H{o} \cite{lp} showed:

\begin{lemma}\mbox{}\\
For $w=+1$ all solutions are the following:\\
For any $t>0,t\neq 3$: $(p,q)=(1,0),(0,1)$.\\
For $t=4$: $(p,q)=(2,3),(3,-2)$.\\
\\
For $w=-1$ all solutions are the following:\\
For $t=1$: $(p,q)=(1,2),(2,-1)$.\\
\\
For $w=+4$ all solutions are the following:\\
For $t=1$: $(p,q)=(3,1),(1,-3)$.\\
\\
For $w=-4$ all solutions are the following:\\
For any $t>0,t\neq 3$: $(p,q)=(1,1),(1,-1)$.\\
For $t=4$: $(p,q)=(5,1),(1,-5)$.\\
\label{thuelemma}
\end{lemma}

Note that congruence consideration mod 8 shows that equation (\ref{thuefamily})
is not solvable for $w=\pm 2$. Further, 
\[
{\rm if} \;\; F_t(p,q)=-4c\;\; {\rm then}\;\; F_t\left(\frac{p-q}{2},\frac{p+q}{2}\right)=c,
\]
(in this case indeed $p$ and $q$ have the same parity), and 
\[
{\rm if} \;\; F_t(p,q)=c\;\; {\rm then}\;\; F_t(p-q,p+q)=4c.
\]
Using Lemma \ref{thuelemma} and the above notes we can figure out the solutions of
(\ref{thuefamily}) for all $w$ being a power of 2 or its negative, and just this
is what we need in our calculation.

\subsection{Index form equations in arbitrary quartic fields}

In this section we detail the method of I.Ga\'{a}l, A.Peth\H{o} and M.Pohst \cite{gppsim}
(see also I.Ga\'{a}l \cite{gaalbook}) which will play an essential role 
in our arguments.

Let $K$ be a quartic field generated by a root
$\xi$ with minimal polynomial 
$f(x)=x^4+a_1x^3+a_2x^2+a_3x+a_4\in\Z[x]$.
We represent any $\alpha\in\Z_K$ in the form
\begin{equation}
\alpha=\frac{1}{d}\left( a+x\xi+y\xi^2+z\xi^3  \right)
\label{alpha8}
\end{equation}
with coefficients 
$a,x,y,z\in\Z$ and with a common denominator $d\in\Z$.
Let $n=I(\xi)$,
\[
F(u,v)=u^3-a_2u^2v+(a_1a_3-4a_4)uv^2+(4a_2a_4-a_3^2-a_1^2a_4)v^3
\]
a binary cubic form over $\Z$ and
\begin{eqnarray*}
Q_1(x,y,z)&=&x^2-xya_1+y^2a_2+xz(a_1^2-2a_2)+yz(a_3-a_1a_2)+z^2(-a_1a_3+a_2^2+a_4)\\
Q_2(x,y,z)&=&y^2-xz-a_1yz+z^2a_2
\end{eqnarray*}
ternary quadratic forms over $\Z$.

\begin{lemma}
If $\alpha$ of (\ref{alpha8}) satisfies
\begin{equation}
I(\alpha)=m,
\label{im}
\end{equation}
then there is a solution $(u,v)\in \Z^2$ of
\begin{equation}
F(u,v)=\pm \frac{d^{6}m}{n}
\label{res}
\end{equation}
such that
\begin{eqnarray}
Q_1(x,y,z)&=&u, \nonumber \\
Q_2(x,y,z)&=&v.   \label{Q12}
\end{eqnarray}
\label{lemma4}
\end{lemma}

In \cite{gppsim} an algorithm is also given for the resolution of the
system of equations (\ref{Q12}) which we shall apply in the following.

\subsection{Index form equations in simplest quartic fields}

Using the coefficients of the polynomial (\ref{pp}), in Lemma \ref{lemma4} we substitute
\[
a_1=-t, a_2=-6, a_3=t, a_4=1,
\]
then we become
\begin{eqnarray}
F(u,v)&=&(u+2v)(u^2+4uv-v^2(t^2+12))\label{ff}\\
Q_1(x,y,z)&=&x^2+txy-6y^2+(t^2+12)xz-5tyz+(t^2+37)z^2\label{qq1}\\
Q_2(x,y,z)&=&y^2-xz+tyz-6z^2\label{qq2}
\end{eqnarray}

We deal with all four cases ($v_2(t)=0,1,2,\geq 3$) in Lemma \ref{kimlemma} parallely.
According to Lemma \ref{lemma4}, for a given $m$, in order to determine
the elements of index $m$ we first have to solve the equation
\[
F(u,v)=\frac{g^6m}{I(\xi)}
\]
with 
\[
F(u,v)=(u+2v)(u^2+4uv-v^2(t^2+12)).
\]
Here $g=2,2,4,4$ and $I(\xi)=2,4,8,16$ according as
$v_2(t)=0,1,2,\geq 3$. We proceed by taking $m=1,2,\ldots,I(\xi)$
until we find solutions.

We write $g^6m/I(\xi)$ in the form $a^\cdot 2^{\ell}$ with
an odd $a\in\{1,3,5,7,9,11,13,15\}$ and $4\leq \ell\leq 12$.
We confer 
\begin{equation}
u+2v=\pm a_1\cdot 2^i
\label{u2v}
\end{equation}
\[
u^2+4uv-v^2(t^2+12)=\pm a_2\cdot 2^{\ell-i}
\]
with odd numbers $a_1,a_2$ satisfying $a_1a_2=a$ and $i=0,1,\ldots ,\ell$.
We become
\begin{equation}
a_1^2\cdot 2^{2i}\pm a_2\cdot 2^{\ell-i}=v^2(t^2+16).
\label{a12}
\end{equation}
The left side $v^2(t^2+16)$ is either zero or
positive.

\noindent
{\bf Case I.}
In case $v^2(t^2+16)=0$ we get $v=0$ for arbitrary $t$.
(\ref{a12}) implies
\[
\frac{a_1^2}{a_2}=2^{\ell -3i}
\]
which is only possible for $a_1=a_2=1$ and $(i,\ell)=(2,6),(3,9),(4,12)$.
By $v=0$ equation (\ref{u2v}) implies $u=\pm 2^i$, therefore (following the
method of \cite{gppsim})
\[
Q_0(x,y,z)=vQ_1(x,y,z)-uQ_2(x,y,z)=0
\]
whence
\[
2^i (y^2-xz+yzt-6z^2)=0.
\]
A non-trivial solution of this quadratic equation is 
$(x_0,y_0,z_0)=(-6,0,1)$. Using an idea of \cite{mordell} Chapter 7
we parametrize all solutions 
$x,y,z$ with rational parameters $p,q,r$ in the form 
\begin{equation}
x=-6p+r,y=q,z=r.
\label{xyz}
\end{equation}
Substituting it into $Q_0(x,y,z)=0$
we obtain 
\[
r(2^ip-2^iqt)=2^iq^2.
\]
Multiplying all equations by $2^ip-2^iqt$
we obtain
\begin{equation}
\begin{array}{cccc}
kx=&2^ip^2&-2^ipqt&-6\cdot 2^iq^2,\\
ky=&      &2^ipq  &-2^itq^2, \\
kz=&      &       &2^iq^2
\end{array}
\label{kxyz}
\end{equation}
with $k\in\Q$. Arrange the coefficients of $p^2,pq,q^2$ on the left hand
sides of the equations above into a 3x3 matrix $C=(c_{ij})$.
By multiplying all equations by the square of the common denominators of
$p,q$ and dividing them by the gcd of the elements of the matrix 
$C=(c_{ij})$ we can
replace $k,p,q$ by integer parameters (cf. \cite{mordell}, \cite{gppsim}) and 
$k$ can be shown to divide 
\[
\frac{\det(C)}{(\gcd(c_{ij}))^2}=2^i
\]
(cf. \cite{gppsim}), \cite{gaalbook}).
Substituting (\ref{kxyz}) into $Q_1(x,y,z)=\pm u$
we obtain
\begin{equation}
2^{2i}\cdot F_t(p,q)=\pm 2^i\cdot k^2
\label{q1}
\end{equation}
with
\[
F_t(p,q)=p^4-tp^3q-6p^2q^2+tpq^3+q^4
\]
where $i=2,3,4$ and $k|2^i$.

This equation is not solvable for $k=1$. 
Using Lemma \ref{thuelemma} and the remarks after it,
solutions exist only for $i=2$ and $i=4$.
The corresponding values are $\ell=6,12,$ respectively.
Using the formula
\[
F(u,v)=\pm\frac{g^6m}{I(\xi)}=2^{\ell}
\]
we can figure out which $v_2(t),I(\xi)$ and $m$ 
may possibly correspond to $i$. The case $m=1$ was deal with
by P.Olajos \cite{olaj}. Using Lemma \ref{thuelemma} we obtain the solutions
$(p,q)$ of equation (\ref{q1}), then $(x,y,z)$ is obtained by (\ref{kxyz}).

\noindent
Finally we obtain the following solutions (and their negatives):\\
For any $t$ with $v_2(t)=0$ or $v_2(t)=1$
\[
(x,y,z)=(2,0,0),(12,2t,-2),(5\pm t,\mp 1+t,-1),(-5\pm t,\mp 1-t,1).
\]
For any $t$ with $v_2(t)=2$ or $v_2(t)\geq 3$
\[
(x,y,z)=(4,0,0),(24,4t,-4),(10\pm 2t,\mp 2+2t,-2), (-10\pm 2t,\mp 2-2t,2).
\]
The corresponding coordinates in the integral basis are listed in the Theorem.

\vspace{1cm}

\noindent
{\bf Case II:}
Using equation (\ref{a12}) and considering 
the possible values of $a_1,a_2,i,\ell$
in case $v^2(t^2+16)>0$ we obtain specific values for $t$.
The possible triples $(t,u,v)$ are listed below.
\par
\noindent
$
\begin{array}{|c|c|c|c|c|c|c|c|c|c|c|c|c|c|c|c|c|}
\hline
t  &2&2&4&4&4&4        &4&4    &4&4&4&4       &7&7&7&7\\ \hline
u  &0&-4&20&-28&16&-32  &12&20  &8&-24&20&-44  &1&3&20&12\\ \hline
v  &1&1&2&2&4&4        &2&2    &4&4&6&6       &-1&-1&-2&2\\ \hline
\end{array}
$
\par
\noindent
$
\begin{array}{|c|c|c|c|c|c|c|c|c|c|c|c|c|c|c|}
\hline
t &8&8&8&8&8&8&8&8&8&8&8&8&8&8\\ \hline
u &10&-14&12&20&22&34&14&-18&-2&-10&34&-54&4&-12 \\ \hline
v &1&1&2&2&3&3&1&1&3&3&5&5&2&2    \\ \hline
\end{array}
$
\par
\noindent
$
\begin{array}{|c|c|c|c|c|c|c|c|c|c|c|c|c|c|c|c|c|}
\hline
t &12&12&12&12   &16&16&16&16&16&16&16&16   &20&20&20&20\\ \hline
u &20&-28&14&-18 &10&-14&28&-36&18&-22&2&-6 &14&-18&36&-44  \\ \hline
v &2&2&1&1       &1&1&2&2&1&1&1&1           &1&1&2&2\\ \hline
\end{array}
$
\par
\noindent
$
\begin{array}{|c|c|c|c|c|c|c|c|c|c|c|c|c|c|c|}
\hline
t &24&24&24&24&24&24  &28&28  &32&32&32&32&32&32   \\ \hline
u &2&-6&44&-52&18&-22 &52&-60 &30&-34&2&-6&60&-68 \\ \hline
v &1&1&2&2&1&1        &2&2    &1&1&1&1&2&2\\ \hline
\end{array}
$
\par
\noindent
$
\begin{array}{|c|c|c|c|c|c|c|c|c|c|c|c|c|c|c|}
\hline
t &48&48&64&64&80&80&96&96     &112&112&128&128&144&144 \\ \hline
u &46&-50&62&-66&78&-82&94&-98 &110&-114&126&-130&142&-146\\ \hline
v &1&1&1&1&1&1&1&1             &1&1&1&1&1&1\\ \hline
\end{array}
$
\par
\noindent
$
\begin{array}{|c|c|c|c|c|}
\hline
t &240&240&256&256 \\ \hline
u &238&-242&254&-258\\ \hline
v &1&1&1&1\\ \hline
\end{array}
$

For each triple $(t,u,v)$ we have to solve the system of equations
(\ref{Q12}). We demonstrate our procedure for $(t,u,v)=(12,20,2)$.
We have
\begin{equation}
Q_1(x,y,z)=x^2+12xy-6y^2+156xz-60yz+181z^2=\pm 20,
\label{ttt1}
\end{equation}
\begin{equation}
Q_2(x,y,z)=y^2-xz+12yz-6z^2=\pm 2,
\label{ttt2}
\end{equation}
\[
Q_0(x,y,z)=2x^2+24xy-32y^2+332xz-360yz+482z^2=0.
\]
The last equation has the non-trivial solution $(x_0,y_0,z_0)=(15,11,-1)$.
We set
\begin{equation}
x=15r+p,\;\;y=11r+q,\;\;z=-r
\label{fxyz}
\end{equation}
(with rational parameters) whence $Q_0(x,y,z)=0$ implies
\[
(8p-16q)r=2p^2+24pq-32q^2.
\]
Multiplying equation (\ref{fxyz}) by $8p-16q$ and using the above equation 
we confer
\begin{eqnarray}
kx&=&-38p^2-344pq+480q^2 \nonumber  \\ 
ky&=&-22p^2-272pq+368q^2 \label{xyzt} \\
kz&=&2p^2+24pq-32q^2\nonumber
\end{eqnarray}
with $k\in\Q$. 
Multiplying all equations by the square of the common denominators of
$p,q$ and dividing them by the gcd of the elements of the above 
coefficient matrix $C=(c_{ij})$ in the right hand side of (\ref{xyzt}) we can
replace $k,p,q$ by integer parameters (cf. \cite{gppsim}).
The number $k$ divides $\det (C)/(\gcd (c_{ij}))^2=3\cdot 2^7$.
Substituting equations (\ref{xyzt}) into the equations
(\ref{ttt1}) and (\ref{ttt2}) we obtain
\[
F_2(p,q)=8p^4+128p^3q+128p^2q^2-3072pq^3+3328q^4=\pm 2\cdot k^2.
\]
We could solve all these equations by using the program package Kash \cite{kash}.
The total CPU time on an average laptop took a couple of hours.
The solutions $(x,y,z)$ are listed below. Note that for $(t,-u,-v)$ we get 
the solutions $(-x,-y,-z)$.
\[
\begin{array}{|c|c|}
\hline
(p,q)&(x,y,z) \\ \hline
(1,0)&(19,11,-1) \\  \hline
(2,1)&(15,11,-1) \\  \hline
(10,-1)&(-5,-37,3) \\  \hline
\end{array}
\]
The corresponding coordinates in the integral basis are listed in the Theorem.

\noindent
Istv\'{a}n Ga\'{a}l and G\'abor Petr\'anyi\\
University of Debrecen, Mathematical Institute \\
            H--4010 Debrecen Pf.12., Hungary \\
            e--mail: igaal@science.unideb.hu, petrarqa@gmail.com  \\ \\


\begin{thebibliography}{10}

\normalsize
\baselineskip=17pt

\bibitem{maple}
B.W.Char, K.O.Geddes, G.H.Gonnet, B.L. Leong, M.B.Monagan, S.M.Watt,
Maple V - language reference manual, New York etc.: Springer-Verlag,
(1991).   Zbl 0758.68038



\bibitem{kash}
M. Daberkow, C.Fieker, J.Kl\"uners, M.Pohst, K.Roegner and K.Wildanger,
KANT V4, J. Symbolic Comput. 24 (1997), 267--283.
http://www.math.tu-berlin.de/$\tilde{\mbox{  } }$kant/.
Zbl 0886.11070, MR1484479

\bibitem{gaalbook}
I.Ga\'al,
Diophantine equations and power integral bases,
Boston, Birkh\"auser, 2002.
Zbl 1016.11059, MR1896601

\bibitem{gaallettl}
I.Ga\'al and G.Lettl, A parametric family of quintic Thue 
equations II., Monatsh. Math., 131 (2000), 29--35.
Zbl 0995.11024, MR1796800

\bibitem{gppsim}
I.Ga\'{a}l, A.Peth\H{o} and M.Pohst, Simultaneous representation
of integers by a pair of ternary quadratic forms -- with an application to
index form equations in quartic number fields, J.Number Theory,
57(1996), 90--104.
Zbl 0853.11023, MR1378574

\bibitem{gpc5}
I.Ga\'al and M.Pohst, Power integral bases in a parametric family of
totally real cyclic quintics, Math. Comp. 66 (1997), 1689-1696.
Zbl 0899.11064, MR1423074

\bibitem{gras}
M.N.Gras, Table num\'erique du nombre de classes et des unit\'es 
des extensions cycliques r\'eelles de degr\'e 4 de $\Q$,
Publ. Math. Fac. Sci. Besançon, Th\'eor. Nombres, Ann\'ee 1977-1978, 
Fasc. 2, 53 p. (1978)
Zbl 0471.12006

\bibitem{jadrievic}
B. Jadrijevi\'c, Establishing the minimal index in a parametric family 
of bicyclic biquadratic fields,
Per. Math. Hungar. 58(2009), 155{-}180.
Zbl 1265.11061, MR2531162

\bibitem{jadrievic2}
B. Jadrijevi\'c, Solving index form equations in the two
parametric families of biquadratic fields, Math. Commun. {\bf 14}(2009),
341--363.
Zbl 1257.11098, MR2743182

\bibitem{kim}
H.K.Kim and J.H.Lee, Evaluation of the Dedekind zeta function at
$s =-1$ of the simplest quartic fields, 
Trends in Math. ,New Ser., Inf. Center for Math. Sci.,
11, No.2 (2009), 63--79.

\bibitem{lazarus}
A.J.Lazarus, On the class number and unit index of simplest quartic fields,
Nagoya Math. J. 121(1991), 1-13.
Zbl 0719.11073, MR1096465


\bibitem{lp}
G.Lettl and A.Peth\H{o},
Complete solution of a family of quartic Thue equations,
Abh. Math. Semin. Univ. Hamb. 65(1995), 365-383.
Zbl 0853.11021, MR1359142

\bibitem{lpv} 
G.Lettl, A.Peth\H{o}, P.Voutier,
Simple families of Thue inequalities,
Trans. Am. Math. Soc. 351(1999), No.5, 1871-1894.
Zbl 0920.11041, MR1487624

\bibitem{mordell}
L.J.Mordell, Diophantine equations,
Pure and Applied Mathematics, 30. London-New York: Academic Press, 1969.
Zbl 0188.34503, MR0249355

\bibitem{olaj}
P.Olajos, Power integral bases in the family of simplest quartic fields,
Exp. Math. 14 (2005), No. 2, 129-132.
Zbl 1092.11042, MR2169516

\bibitem{shanks}
D.Shanks, The simplest cubic fields,
Math. Comput. 28 (1974), 1137-1152.
Zbl 0307.12005, MR0352049



\end{thebibliography}
\end{document}